\documentclass[reqno]{amsart}
\usepackage{amssymb,amsthm,amsfonts,amstext}

\usepackage{amsmath}
\usepackage{newcent}       
\usepackage{helvet}         
\usepackage{courier}        

\makeatletter \@addtoreset{equation}{section} \makeatother

\renewcommand\thefigure{\thesection.\@arabic\c@figure}
\renewcommand\thetable{\thesection.\@arabic\c@table}
\newtheorem{theorem}{Theorem}[section]

\newtheorem{proposition}[theorem]{Proposition}

\newcommand{\mc}[1]{{\mathcal #1}}

\newcommand{\bb}[1]{{\mathbb #1}}
\newcommand{\<}{\langle}
\renewcommand{\>}{\rangle}

\begin{document}
\author{Patr\'{\i}cia Gon\c{c}alves}
\address{CMAT, Centro de Matem\'atica da Universidade do Minho, Campus
de Gualtar, 4710-057 Braga, Portugal} \email{patg@impa.br and patg@math.uminho.pt}

\author{Milton Jara}

\address{{IMPA, Instituto Nacional de Matem\'atica Pura e Aplicada\\
Estrada Dona Castorina 110\\
Jardim Bot\^anico\\
22460-320 Rio de Janeiro-RJ\\
Brazil} and {Ceremade, UMR CNRS 7534,
    Universit\'e Paris - Dauphine,
    Place du Mar\'echal De Lattre De Tassigny
    75775 Paris Cedex 16 - France}
\newline
e-mail: \rm \texttt{mjara@impa.br}}

\title[Crossover to the KPZ equation]{Crossover to the KPZ equation}
\date{\today}


\begin{abstract}
 We characterize the crossover regime to the KPZ equation for a class
of one-dimensional weakly asymmetric exclusion processes. The crossover depends on the strength asymmetry $an^{2-\gamma}$ ($a,\gamma>0$) and it occurs at
$\gamma=1/2$. We show that the density field is a solution of an Ornstein-Uhlenbeck equation if $\gamma\in(1/2,1]$, while for $\gamma=1/2$ it is an energy
solution of the KPZ equation. The corresponding crossover for the current of particles is readily obtained.
\end{abstract}

\maketitle
\section{\small{Introduction}}
 One-dimensional weakly asymmetric exclusion processes arise as simple models for the
growing of random interfaces. For those processes the microscopic dynamics is given by stochastic lattice gases with hard core exclusion with a weak asymmetry to the right. The presence of a weak asymmetry breaks down the detailed balance condition, which forces the system to exhibit a
non trivial behavior even in the stationary situation. Using renormalization group techniques, the dynamical scaling exponent has been
established as $z=3/2$ and one of the challenging problems is to derive the limit distribution of the density and the current of particles \cite{HS}.

 For asymmetric exclusion processes, partial answers have been given in
particular settings as starting the system from the stationary state and from specific initial conditions. For these models, under a certain spatial shifting and time speeding, the current of particles has Tracy-Widom distribution, see \cite{FS,J,TW1,TW2}. The Tracy-Widom
distribution was initially obtained in the context of large N statistics of the largest eigenvalue of random matrices, but has been recently
obtained as the scaling limit of stochastic fields of random models, see \cite{ACQ,BQS,J,SS1,SS2} and references therein.

Here we are interested in establishing the equilibrium density fluctuations for weakly asymmetric exclusion processes with strength asymmetry
$an^{2-\gamma}$, that we fully describe below. We consider the system under the invariant state: a Bernoulli product measure of parameter $\rho\in{[0,1]}$ that we denote by $\nu_\rho$. By relating the current with the density of particles, as a consequence of last result we derive the equilibrium
fluctuations of the current of particles.

 The weakly
asymmetric simple exclusion process was studied in \cite{DMPS,DG}, for strength asymmetry $n$ (that corresponds to $\gamma=1$ in our case), and in \cite{BG} for the strength
asymmetry $n^{3/2}$ (that corresponds to $\gamma=1/2$ in our case). For $\gamma=1$, the equilibrium density fluctuations are given by an Ornstein-Uhlenbeck process which
implies the current to have Gaussian distributions, see \cite{DMPS,DG}. For $\gamma=1/2$, \cite{BG} used the
Cole-Hopf transformation to derive the non-equilibrium fluctuations of the current. By using the Cole-Hopf transformation initially introduced
in \cite{G}, one obtains an exponential process. For this exponential process, the limiting fluctuations are given by the stochastic heat equation, which is a linear equation, making the asymptotic analysis much easier. As a consequence, the fluctuations of the original process are easily
recovered. In our approach, we study the weakly asymmetric exclusion directly, but in order to identify the limiting density field as a weak solution of a stochastic partial differential equation, we have to overcome the difficulty of closing the equation
by means of the Boltzmann-Gibbs principle. For this reason our results are restricted to the equilibrium setting.

 It is known from \cite{DMPS,DG} that for $\gamma=1$, the limit density field is a solution of an Ornstein-Uhlenbeck equation which has a drift term. This drift term comes from the asymmetric part of
the dynamics and can be removed by taking the process moving in a reference frame with constant velocity. By removing the drift of the system,
there is no effect of the strength of the asymmetry on the distribution of the limit density field. In order to see how far this picture can go, we strengthen the asymmetry by decreasing the value of $\gamma$. We can show that, for $\gamma\in{(1/2,1]}$ there is still no effect of the strength of the asymmetry on the limiting density field. In this case, the limiting density field is still solution of the Ornstein-Uhlenbeck equation as for $\gamma=1$ and for this reason the process belongs to the Edwards-Wilkinson \cite{EW} universality class.  Nevertheless, for $\gamma=1/2$ the limiting distribution "feels" the effect of the strengthening of the asymmetry, by developing a non linear term in the
equation that characterizes the limiting density field. In this case the limiting density field is a solution of the Kardar-Parisi-Zhang (KPZ)
equation, so that for $\gamma=1/2$ the process belongs to the KPZ \cite{KPZ} universality class.

The KPZ equation was proposed in \cite{KPZ} to model the growth of random interfaces. Denoting by $h_{t}$ the height of the interface, this
equation reads as
\begin{equation*}
\partial_{t}h=D\Delta h+a(\nabla h)^2+\sigma \mathcal{W}_{t},
\end{equation*}
 where $D, a, \sigma$ are related to the thermodynamical
properties of the interface and $\mathcal{W}_{t}$ is a Gaussian space-time white noise with covariance given by
\begin{equation*}
E[\mathcal{W}_{t}(u)\mathcal{W}_{s}(v)]=\delta(t-s)\delta(u-v).
\end{equation*}

 According to the dynamical scaling exponent $z=3/2$, a non-trivial behavior occurs under the scaling $h_{n}(t,x)=n^{-1/2}h(tn^{3/2},x/n)$.
 This means, roughly speaking, that in our case, for $\gamma=1/2$ a non trivial behavior is expected even in the stationary situation and in that case the model belongs
to the universality class of the KPZ equation.

To our knowledge, a rigorous mathematical proof of the characterization of the intermediate state
between the Ornstein-Uhlenbeck process and the crossover to the KPZ equation, was lacking so far, see \cite{PM} and references therein. Nevertheless we refer the reader to the paper \cite{ACQ} in which the authors characterize the crossover regime for a special case of weakly asymmetric exclusion process, but going through the Cole-Hopf transformation.

As a consequence of obtaining the fluctuations of the density of particles, we obtain the equilibrium fluctuations for the current of particles for different strength regimes depending on the strength asymmetry. More precisely, we show that for $\gamma>1/2$ the current properly centered and re-scaled converges to a fractional Brownian motion with Hurst parameter $1/4$ and for $\gamma=1/2$ the limit process is given in terms of the solution of the KPZ equation.

The existence of a non-trivial crossover regime for weakly asymmetric systems was found in \cite{BD}. There, the authors use the theory developed
in \cite{BDSGJLL1} and \cite{BDSGJLL2}, and found that for a phase of weak asymmetry the fluctuations of the current are Gaussian, while in
the presence of stronger asymmetry they become non-Gaussian. Here, we provide the characterization of the transition from the Edwards-Wilkinson class to the KPZ class, for general weakly asymmetric exclusion processes. We prove that the transition depends on the
strength of the asymmetry without having any other intermediate state and by establishing precisely the strength in order to have the crossover. We point out here, that our results are also valid for weakly asymmetric exclusion processes with finite-range interactions. All the proofs follow in this case with minor notational modifications.

Here follows an outline of this paper. In the second section, we introduce the model and we describe the equilibrium density and current
fluctuations for the process under different strength asymmetry regimes. In the third section, we sketch the proof of the results for the
intermediate state regime and in the fourth section we recall briefly the results about the crossover to the KPZ equation established in \cite{GJ5}.

\section{\small{Equilibrium fluctuations}}

Let $\eta_t$ be the weakly asymmetric exclusion process evolving
on $\mathbb{Z}$. The state space of this Markov process is $\Omega:=\{0,1\}^{\mathbb{Z}}$ and its dynamics can be described as follows. On a configuration $\eta\in{\Omega}$ and
after a mean one exponential time, a particle jumps to an empty
neighboring site according to a transition rate that has a weak asymmetry to the right and depends on a function $c(\eta)$. We assume
$c:\Omega\rightarrow{\mathbb{R}}$ to be a {\em{local}} function, {\em{bounded from above and below}}, and
that turns the system {\em{gradient}} and {\em{reversible}} with respect to the
stationary state $\nu_\rho$. The gradient condition is the most restrictive one
and requires the existence of a local function $h:\Omega\rightarrow{\mathbb{R}}$ such that for any $\eta \in \Omega$
\begin{equation*}
c(\eta)(\eta(1)-\eta(0)) = \tau_1 h(\eta)- h(\eta).
\end{equation*}

 The process is speeded up on the diffusive time scale $n^2$ so that
$\eta_t^n=\eta_{tn^2}$.
Here and in the sequel, for $\eta\in{\Omega}$ and $x\in{\mathbb{Z}}$ we denote by $\tau_x \eta$ the space translation by $x$, namely for
$y\in{\mathbb{Z}}$, $\tau_x\eta(y)=\eta(y+x)$ and for a function $f:\Omega \to \mathbb{R}$ we denote by $\tau_x f(\eta)$ the induced translation
in $f$, namely $\tau_x f(\eta)=f(\tau_x \eta)$. Let $c_x(\eta) = \tau_x c(\eta)$. For the configuration $\eta$, the transition rate from $x$ to $x+1$ is given by $c_{x}(\eta)p_n$  and from $x+1$ to
$x$ is given by $c_{x}(\eta)q_n$, where $p_n:=(1+a/n^\gamma)/2$, $q_n:=1-p_n$ and $a>0$. We refer the reader to \cite{GJ5} for a complete discussion about the assumptions on
$c(\cdot)$.

 The generator of this process acts over local functions $f:\Omega \to \mathbb{R}$ as
\[
\mathcal{L}_n f(\eta) = n^2 \sum_{x \in \bb Z} c_x(\eta) \big\{p_n \eta(x)(1-\eta(x+1)) + q_n \eta(x+1)(1-\eta(x))\big\} \nabla_{x,x+1} f(\eta),
\]
where $n \in \bb N$, $\nabla_{x,x+1} f(\eta) = f(\eta^{x,x+1})-f(\eta)$ and
\[
\eta^{x,x+1}(z) =
\begin{cases}
\eta(x+1),& z=x\\
\eta(x), & z=x+1\\
\eta(z), & z\neq x,x+1.\\
\end{cases}
\]

If $a=0$ the process $\eta^{n}_{t}$ is said to be symmetric and if $c(\cdot) \equiv 1$ the process is the symmetric simple
exclusion process. On the other hand if $c(\cdot) \equiv 1$, $a=1$ and $\gamma=1$, the process is the weakly asymmetric simple exclusion
process, studied in \cite{DMPS,DG}, and that can be interpreted as having the symmetric and asymmetric dynamics speeded up by $n^2$ and $n$, respectively.

We notice that decreasing the value of $\gamma$ in $p_n$ above, corresponds
to speeding up the asymmetric part of the dynamics on longer time scales as $n^{2-\gamma}$, while the scaling of the symmetric dynamics is not affected.

 A stationary
state for this process is the Bernoulli product measure on $\Omega$ of parameter $\rho\in [0,1]$ that we denote by $\nu_{\rho}$ and whose marginal at $\eta(x)$ is given by $\nu_{\rho}(\eta:\eta(x)=1)=\rho.$

\subsection{Hydrodynamic Limit}

Here we recall briefly the {\em{hydrodynamic limit}} for
$\eta_{t}^{n}$. The hydrodynamical scaling corresponds to the strength asymmetry $n$.

For that purpose we introduce the empirical measure as the positive
measure in $\bb R$ defined by
\begin{equation*}
\pi_t^n(dx) = \frac{1}{n} \sum_{x \in \mathbb{Z}} \eta_t^n(x) \delta_{x/n}(dx),
\end{equation*}
 where for $u\in{\bb R}$, $\delta_{u}$ is the Dirac measure at
$u$.

Take $\rho_{0}:\mathbb{R}\rightarrow{[0,1]}$ a strictly positive and piecewise continuous function such that there exists $\rho \in (0,1)$ satisfying $\int |\rho_0(x)
-\rho|dx<+\infty$. Start the process $\eta_t^n$ from $\{\mu_n;n \in \bb N\}$ - a product measure on $\Omega$, whose marginal at $\eta(x)$ is Bernoulli of parameter
$\rho_{0}(x/n)$, namely:
\begin{equation*}
\mu_{n}(\eta:\eta(x)=1)=\rho_0(x/n).
\end{equation*}

 Then, $\pi_t^n(dx)$ converges in probability to the deterministic measure $\rho(t,x) dx$, where $\{\rho(t,x); t \geq 0, x \in
\bb R\}$ is the unique weak solution of the {\em{viscous Burgers equation}}
\begin{equation*}
\left\{
\begin{array}{l}
\partial_t \rho(t,x) = \frac{1}{2} \Delta \varphi_h(\rho(t,x)) -a\nabla
\beta(\rho(t,x)) \\
\quad\\
\rho(0,u) = \rho_0(u), u\in{\mathbb{R}}.
\end{array}
\right.
\end{equation*}
 where $\varphi_h(\rho) = \int h d\nu_\rho$, $\beta(\rho) = \chi(\rho)\int c d\nu_\rho$ and $\chi(\rho)=\rho(1-\rho)$. We notice that, thanks to the gradient
(reversible) condition of $c(\cdot)$, $\varphi_h(\rho)$ (respectively $a\beta(\rho)$) corresponds to the expectation with respect to $\nu_{\rho}$ of the
instantaneous current of the symmetric (respectively asymmetric) part of the dynamics. Our purpose here is to analyze the fluctuations of the empirical
measure from the stationary state $\nu_{\rho}$ and from there, to derive the fluctuations of the current of particles.

\subsection{Equilibrium Fluctuations}

From now on we fix a density $\rho \in (0,1)$. Denote by $\mathcal{S}(\bb R)$ the Schwartz space, i.e. the space of rapidly decreasing functions and denote by $\mathcal{S}'(\bb R)$ its dual with respect to the inner product of $L^2(\mathbb{R})$. Let $\mc D([0,\infty), \mc S'(\bb R))$ be the space of $c\grave{a}dl\grave{a}g$ trajectories from $[0,\infty)$ to $\mathcal{S}'(\bb R)$.  Define $\{ \mc Y_t^{n}; t \geq 0\}$ as the {\em{density fluctuation field}}, a linear functional acting on $H\in{\mathcal{S}(\bb R)}$ as
\begin{equation}\label{density field}
\mc Y_t^{n}(H) = \frac{1}{\sqrt n} \sum_{x \in \bb Z}
H\Big(\frac{x}{n}\Big) \big(\eta_t^n(x)-\rho \big).
\end{equation}
We start by considering $\gamma=1$, which as mentioned above corresponds to the hydrodynamic strength asymmetry. By computing the
characteristic function of $\mc Y_0^{n}(H)$ it follows that $\mc Y_0^{n}$ converges in distribution to a spatial white noise of variance
$\chi(\rho)$ - {the  \em{static compressibility of the system}}. Moreover, $\{\mc Y_t^{n};n \in \bb N\}$ converges in distribution with respect to
the {\em{Skorohod topology}} of the space of $c\grave{a}dl\grave{a}g$ trajectories $\mc D([0,\infty), \mc S'(\bb R))$ to the process $\mc Y_t$, solution of the Ornstein-Uhlenbeck equation
\begin{equation}
\label{OU1}
d \mc Y_t = \frac{\varphi_h'(\rho)}{2}  \Delta \mc Y_t dt -a\beta'(\rho)
\nabla \mc Y_t dt
       +\sqrt{\beta(\rho)} \nabla d \mc W_t,
\end{equation}
where $\mc W_t$ is a space-time white noise of unit variance. To prove this result, since the behavior at $t=0$ is characterized it remains to
analyze the time evolution of the limit density field. For this purpose we introduce the martingale associated to \eqref{density field} and we analyze
its asymptotic behavior. The martingale decomposition presents two integral terms which cannot be written as a function of the density
fluctuation field given in \eqref{density field}. As a consequence we cannot identify straightforwardly the limit density field as a weak solution to some
stochastic partial differential equation. In order to perform this identification, a replacement argument is required. This replacement is known
as the {\em Boltzmann-Gibbs} principle and was introduced in \cite{Ros}. With this result in hand and since the system is gradient, the equilibrium fluctuations are easily
derived, see for example \cite{DMPS} for a detailed proof.

 Now we analyze equation (\ref{OU1}). This equation has a drift term
that vanishes if $\beta'(\rho)$ is equal to $0$, and in that case the
equation does not depend on $a$ any longer. In order to remove the
drift term (which arises from the asymmetric part of the dynamics)
from the limiting density field, we take $\eta_{t}^{n}$ moving in a reference
frame with constant velocity $a\beta'(\rho)n$. In order to see the
dependence on the strength of the asymmetry on the limit process, we redefine the
{\em{density fluctuation field}} on $H\in\mathcal{S}(\mathbb{R})$ as:
\begin{equation}\label{gamma density field}
\mc Y_t^{n,\gamma}(H) = \frac{1}{\sqrt n} \sum_{x \in \bb Z}
T_{t}^{\gamma}H\Big(\frac{x}{n}\Big) \big(\eta_t^n(x) -\rho\big),
\end{equation}
where $T_{t}^{\gamma}H(\cdot)=H(\cdot-a\beta'(\rho)tn^{1-\gamma})$ and $a \beta'(\rho)n^{2-\gamma}$ is the velocity of the system (be aware of
the scaling parameter $x\rightarrow{x/n}$). As above, it is not hard to show that for $\gamma=1$, $\{\mc Y_t^{n,\gamma}; n\in \bb N\}$ converges in the Skorohod topology of $\mc D([0,\infty), \mc S'(\bb R))$ to the process $\mc
Y_t$ solution of the Ornstein-Uhlenbeck equation
\begin{equation}
\label{OU2} d\mc Y_t = \frac{ \varphi_h'(\rho)}{2} \Delta \mc Y_t dt + \sqrt{\beta(\rho)} \nabla d \mc W_t,
\end{equation}
which corresponds to equation \eqref{OU1} with $a=0$. Since for $\gamma=1$ the fluctuations of the density are given by \eqref{OU2}, the system belongs to the Edwards-Wilkinson universality class \cite{EW}.

In order to see the effect of incrementing the strength asymmetry in the limit density field we decrease the value of
$\gamma$. As discussed in \cite{BG}, the effect of the asymmetry is presented in the limiting density field when $\gamma=1/2$ and in that case $\mc
Y_t$ has a very different qualitatively behavior from the one obtained for $\gamma=1$, namely the solution of \eqref{OU2}.
Here we
characterize the limiting density field $\mathcal{Y}_{t}$ for the intermediate state regime, namely we show that:

\begin{theorem} \label{flu}
If $\gamma\in{(1/2,1]}$, then the sequence $\{\mc Y_t^{n,\gamma};n \in \bb N\}$ converges in distribution with respect to Skorohod topology
of $\mc D([0,\infty), \mc S'(\bb R))$ to the process $\mc Y_t$ solution of the Ornstein-Uhlenbeck equation \eqref{OU2}.
\end{theorem}

By the previous result, we have that for $\gamma\in(1/2,1]$ the system still belongs to the Edwards-Wilkinson universality class. So, if we want to see the
effect of strengthening the asymmetry in the limiting density field, we have to take $\gamma=1/2$ which is in agreement with the result of \cite{BG}. Recently in \cite{GJ5}, it was shown that for $\gamma=1/2$, the sequence $\{\mc Y_t^{n,\gamma};n \in \bb N\}$ is tight and any limit point
is an energy solution of the KPZ equation:
\begin{equation}\label{sB}
d\mathcal{Y}_{t}=\frac{\varphi_h'(\rho)}{2}\Delta \mathcal{Y}_{t}dt-\frac{a \beta''(\rho)}{2}\nabla( \mathcal{Y}_{t})^2dt+\sqrt{\beta(\rho)}\nabla
d\mathcal{W}_{t}.
\end{equation}
So, our result says that weakly asymmetric simple exclusion processes belong to the Edwards-Wilkinson universality class for $\gamma\in(1/2,1]$; and cross to
the KPZ universality class when the strength asymmetry is precisely given by $an^{3/2}$, see \cite{GJ5}.

\section{\small{Beyond the hydrodynamic time scale}}

In this section we give an outline of the proof of the equilibrium density and current fluctuations in the intermediate state regime (with $\gamma\in{(1/2,1]}$),
i.e. between the Ornstein-Uhlenbeck process \eqref{OU2} and the KPZ equation \eqref{sB}.

\subsection{Density fluctuations}
Here we prove Theorem \ref{flu}. Recall the definition of the density fluctuation field given in \eqref{gamma density field}. In this setting,
we remove the drift of the process, so that we suppose to have $\eta_{t}^{n}$ moving in a reference frame with constant velocity given by
$a\beta'(\rho)n^{2-\gamma}$ with $\gamma\in(1/2,1]$.

At $t=0$, by computing the characteristic function of $\mathcal{Y}_{0}^{n,\gamma}$, this field converges to a spatial white noise of variance
$\chi(\rho)$. Now, we analyze the asymptotic behavior of some martingales associated to $\mathcal{Y}_{t}^{n,\gamma}$, in order to identify the limit density field
 $\mathcal{Y}_{t}$ as a weak solution of the stochastic partial differential equation \eqref{OU2}.
For that purpose, fix $H\in{\mathcal{S}(\bb R)}$ and notice that by Dynkin's formula,

\begin{equation*}
\mathcal{M}_t^{n,\gamma}(H)=\mc Y_t^{n,\gamma}(H)-\mc Y_0^{n,\gamma}(H) -\mathcal{I}^{n,\gamma}_{t}(H)-\mathcal{A}_{t}^{n,\gamma}(H)
\end{equation*}
 is a martingale with
respect to the natural filtration $\mathcal{F}_{t}=\sigma(\eta_s,s\leq{t})$, where
\begin{equation*}
\begin{split}
&\mathcal{I}_{t}^{n,\gamma}(H)=\int_0^t \frac{1}{2\sqrt n} \sum_{x \in \bb Z} \Delta^nT_{s}^{\gamma}H\Big(\frac{x}{n}\Big) \Big(\tau_x
h(\eta_s^n)-\varphi_h(\rho)\Big) ds,\\
& \mathcal{A}_{t}^{n,\gamma}(H)= \int_0^t \frac{n^{1-\gamma}}{\sqrt n} \sum_{x \in \bb Z} \nabla^n
T_{s}^{\gamma}H\Big(\frac{x}{n}\Big)\tau_{x}V_{f}(\eta_{s}^{n})ds,\\
 &f(\eta)=ac(\eta)(\eta(1)-\eta(0))^2/2,\\
&\tau_{x}V_{f}(\eta)=\tau_x f(\eta)-a\beta(\rho)-a\beta'(\rho)(\eta(x)-\rho),\\
 & \Delta^n H(x/n):=n^2(H((x+1)/n)+H((x-1)/n)-2H(x/n)) \hspace{0,5cm} and\\
 &\nabla^nH(x):=n(H((x+1)/n)-H(x/n)).
\end{split}
\end{equation*}
Notice that $\Delta^n$ and $\nabla^n$ are the discrete Laplacian and the discrete derivative, respectively. We point out that the mean of $f$
with respect to $\nu_{\rho}$ is given by $a\beta(\rho)$.

The quadratic variation of $\mathcal{M}_t^{n,\gamma}(H)$ equals to
\begin{equation*}
\<\mathcal{M}^{n,\gamma}(H)\>_t=\int_0^t \frac{1}{2n} \sum_{x \in \bb Z} \Big(\nabla^n
T_{s}^{\gamma}H\Big(\frac{x}{n}\Big)\Big)^2\tau_xf(\eta_{s}^{n})\Big(1+\frac{a}{n^\gamma}\Big)ds.
\end{equation*}

We notice that if $\gamma>0$, the term corresponding to $a/n^\gamma$ inside last integral, vanishes in $L^2(\mathbb{P}_{\nu_{\rho}})$ as $n\rightarrow{+\infty}$. From now on, $\mathbb{P}_{\nu_{\rho}}$ denotes the distribution of the Markov process $\eta_{t}^n$
 starting from the stationary state $\nu_{\rho}$ and $\mathbb{E}_{\nu_{\rho}}$ denotes the expectation with respect to $\mathbb{P}_{\nu_{\rho}}$. Last term arises
from the asymmetric part of the dynamics and for this reason, only for a big strength of the asymmetry (namely for $\gamma=0$; which corresponds
to speeding up the asymmetric part of the dynamics by $n^2$) it will give rise to additional stochastic fluctuations of the system. From
simple computations, it follows that the limit as $n\rightarrow{+\infty}$ of the martingale $\mathcal{M}_{t}^{n,\gamma}(H)$ is given by $||\sqrt{\beta(\rho)}\nabla
H||_{2} \mathcal{W}_{t}(H)$, where $\mathcal{W}_{t}(H)$ is a Brownian motion and $||\cdot||_{2}$ denotes the $L^{2}(\mathbb{R})$-norm.

Now, we need to analyze the limit of the integral terms. We start by the less demanding, namely $\mathcal{I}_{t}^{n,\gamma}(H)$. Invoking the
Boltzmann-Gibbs principle introduced in \cite{Ros}, $\mathcal{I}_{t}^{n,\gamma}(H)$ can be written as
\begin{equation*}
\int_0^t \mathcal{Y}_{s}^{n,\gamma}\Big(\frac{ \varphi_h'(\rho)}{2}\Delta^n H\Big) ds
\end{equation*}
 plus an $L^2(\bb P_{\nu_{\rho}})$ negligible term.

Now we analyze $\mathcal{A}_{t}^{n,\gamma}(H)$. We notice that if we were considering the case $\gamma=1$ the Boltzmann-Gibbs principle as stated
in \cite{Ros}, would be saying that $\mathcal{A}_{t}^{n,\gamma}(H)$ vanishes in $L^2(\mathbb{P}_{\nu_{\rho}})$ as $n\rightarrow{+\infty}$. Since $\gamma<1$, the result in \cite{Ros}
does not give us the information we need about the limit of this integral term. Nevertheless, according to Corollary 7.4 of \cite{G1}, it follows that
in fact for $\gamma\in{(1/2,1]}$ the same result is true. Indeed the integral term $\mathcal{A}_{t}^{n,\gamma}(H)$ still vanishes in $L^2(\mathbb{P}_{\nu_{\rho}})$ as
$n\rightarrow{+\infty}$ for $\gamma\in{(1/2,1]}$. We remark that the mentioned result in \cite{G1} was proved for the symmetric simple exclusion
process but it is also true for our model. More explicitly, for our case that result says the following:

\begin{proposition}[Stronger Boltzmann-Gibbs Principle \cite{G1}]\label{sBG}
\quad\

 Let $\psi:\Omega \to \bb R$ be a local function and let $\varphi_\psi(\rho):=E_{\nu_{\rho}}[\psi(\eta)]$. For $\gamma\in{(1/2,1]}$ and $H \in \mc S(\bb R)$,
 it holds that
\begin{equation*}
\lim_{n \to \infty} \bb E_{\nu_{\rho}}\Big[\Big(\int_0^t
\frac{n^{1-\gamma}}{\sqrt n} \sum_{x \in \bb Z}
H\Big(\frac{x}{n}\Big)\big(\tau_x \psi(\eta_s^n)-\varphi_\psi(\rho)-\varphi_\psi'(\rho)(\eta_s^n(x)-\rho)\big)ds\Big)^2\Big]=0.
\end{equation*}
\end{proposition}

 In \cite{G1}, this result was proved for $\psi(\eta)=(\eta(0)-\rho)(\eta(1)-\rho)$, but the proof holds for any local function, since the fundamental ingredients invoked along the proof, namely the spectral gap bound and the equivalence of
ensembles, hold in general for local functions. As a consequence of last result, $\mathcal{A}_{t}^{n,\gamma}(H)$ still vanishes in $L^2(\mathbb{P}_{\nu_{\rho}})$ as $n\rightarrow{+\infty}$ for
$\gamma\in{(1/2,1]}$.

Putting together the previous observations, for $\gamma\in{(1/2,1]}$, the limit density field $\mathcal{Y}_{t}(H)$ satisfies:
\begin{equation*}
\mathcal{Y}_{t}(H)=\mc Y_0(H)+ \int_0^t \mathcal{Y}_{s}\Big(\frac{ \varphi_h'(\rho)}{2}\Delta H\Big) ds+||\sqrt{\beta(\rho)}\nabla H||_{2}
\mathcal{W}_{t}(H),
\end{equation*}
so that $\mathcal{Y}_{t}$ is a weak solution of (\ref{OU2}).
With this decomposition it follows that the covariance of the limit field is given on $H,G\in{\mathcal{S}(\mathbb{R})}$ by
\begin{equation}\label{covariance}
E[\mathcal{Y}_t(H)\mathcal{Y}_s(G)]=\chi(\rho)\int_{\mathbb{R}}T_{t-s}H(x)G(x)dx,
\end{equation}
where $\{T_{t}\}_{t\geq{0}}$ is the semigroup associated to the operator $\frac{ \varphi_h'(\rho)}{2}\Delta$.

In order to complete the argument, it remains to show tightness of the sequence
$\{\mathcal{Y}_{t}^{n,\gamma}\}_{n\in{\mathbb{N}}}$. In \cite{G1} this was shown for the asymmetric simple exclusion process, but the same computations with minor modifications, hold for the
processes we consider here.

\subsection{Current fluctuations}

Now we want to derive the fluctuations of the current of particles, from the fluctuations of the density of particles. Since we took the process
moving in a reference frame with constant velocity $a\beta'(\rho)n^{2-\gamma}$, we consider the current of particles through a moving bond.

Fix $\gamma\in{(1/2,1]}$. For a site $x$ denote by $\mathcal{J}_{x}^{n}(t)$, the current of particles through the bond $\{x,x+1\}$, i.e. the
number of particles that jump from the site $x$ to $x+1$, minus the number of particles that jump from $x+1$ to $x$ during the time
interval $[0,tn^2]$. Formally
\begin{equation*}
\mathcal{J}_{x}^{n}(t)=\sum_{y\geq{x+1}}\Big(\eta^{n}_{t}(y)-\eta^{n}_{0}(y)\Big),
\end{equation*}
 so that it is the difference between the density fluctuation field \eqref{density
field} at time $t$ and at time $0$, evaluated on the Heaviside function $H_{x}=1_{(x,\infty)}$.

Consider the line $a_{x}=x+[a\beta'(\rho)tn^{1-\gamma}]$ and let $\mathcal{J}_{x}^{n,\gamma}(t)$ be the current of particles through the
time-dependent bond $\{a_{x},a_{x}+1\}$. For $x\in{\bb R}$, $[x]$ denotes the integer part of $x$. Take $x=0$ to simplify the exposition but for
any other site the results stated below are also true.

Up to the strength asymmetry $n^{2-\gamma}$ with $\gamma>1/2$, we are able to show that the current properly centered and re-scaled converges to a fractional Brownian motion of Hurst parameter $H=1/4$:

\begin{theorem}
Fix $x\in{\mathbb{Z}}$, $\gamma\in{(1/2,1]}$ and let
\begin{equation*}
\mathcal{Z}_{t}^{n}=\frac{1}{\sqrt{n}}\Big\{\mathcal{J}_{x}^{n,\gamma}(t)-\mathbb{E}_{\nu_{\rho}}[\mathcal{J}_{x}^{n,\gamma}(t)]\Big\}.
\end{equation*}
 Then, for
every $k\geq{1}$ and every $0\leq{t_{1}}<{t_{2}}<...<t_{k}$, $(\mathcal{Z}_{t_{1}}^{n},...,\mathcal{Z}_{t_{k}}^{n})$ converges in law to a Gaussian vector
$(\mathcal{Z}_{t_{1}},...,\mathcal{Z}_{t_{k}})$ with mean zero and covariance given by
\begin{equation*}
E[\mathcal{Z}_{t}\mathcal{Z}_{s}]=\sqrt{\frac{2\varphi_h'(\rho)}{\pi}}\chi(\rho)(\sqrt{t}+\sqrt{s}-\sqrt{t-s})
\end{equation*}
 provided $s\leq{t}$.

\end{theorem}

 We sketch here the proof of last result. The idea of the argument is to obtain the fluctuations of the current from the fluctuations of the
density of particles, namely from Theorem \ref{flu}. This argument was initially proposed in \cite{RV}. For more details we refer the reader to
\cite{G1} and \cite{JL}.

From the stronger Boltzmann-Gibbs principle as stated in Proposition \ref{sBG}, it follows that for $\gamma\in{(1/2,1]}$
\begin{equation*}
\mathcal{J}_{0}^{n,\gamma}(t)-(\mathcal{Y}_{t}^{n,\gamma}(G_{\ell})-\mathcal{Y}_{0}^{n,\gamma}(G_{\ell})),
\end{equation*}
 vanishes in $L^2(\mathbb{P}_{\nu_{\rho}})$
as $\ell\rightarrow{+\infty}$ uniformly over $n$, where $\{G_{\ell}\}_{\ell\in{\mathbb{N}}}$ is a sequence approximating the Heaviside function $H_{0}$, which can be taken for example equal to
$G_{\ell}(x)=(1-x/\ell)^+$.

Combining last result with the convergence of $\mathcal{Y}_{t}^{n,\gamma}$, it follows that $\mathcal{Z}_{t}^{n}$
converges to a random variable, which
formally reads as
\begin{equation*}
\mathcal{Y}_{t}(H_{0})-\mathcal{Y}_{0}(H_{0}),
\end{equation*}
where $\mathcal{Y}_t$ is the solution of the Ornstein-Uhlenbeck equation \eqref{OU2}.

The same argument can be applied to show the same result for any vector $(\mathcal{Z}_{t_{1}},...,\mathcal{Z}_{t_{k}})$. To compute the covariance we do the following:
\begin{equation*}
\begin{split}
E[\mathcal{Z}_{t}\mathcal{Z}_s]&=E[\{\mathcal{Y}_{t}(H_{0})-\mathcal{Y}_{0}(H_{0})\}\{\mathcal{Y}_{s}(H_{0})-\mathcal{Y}_{0}(H_{0})\}]\\
&=\lim_{\ell\rightarrow{+\infty}}E[\{\mathcal{Y}_{t}(G_{\ell})-\mathcal{Y}_{0}(G_{\ell})\}\{\mathcal{Y}_{s}(G_{\ell})-\mathcal{Y}_{0}(G_{\ell})\}].
\end{split}
\end{equation*}
Now we use \eqref{covariance} to write last expression as
\begin{equation*}
\chi(\rho)\lim_{l\rightarrow{+\infty}}\int_{\mathbb{R}}\Big(T_{t-s}G_{\ell}(x)G_{\ell}(x)-G_{\ell}(x)T_tG_{\ell}(x)-G_{\ell}(x)T_sG_{\ell}(x)+G_{\ell}^2(x)\Big)dx.
\end{equation*}
Using the definition of $T_{t}(G_{\ell})$ we get the covariance stated in the theorem.

We notice that last convergence takes place in the sense of finite-dimensional distributions. Since the distributions of $\mathcal{Y}_{t}(H_{0})$ are Gaussian, this implies the limit current to be Gaussian distributed.

\section{The Crossover regime}

Here we describe briefly the limit density fluctuation field for $\gamma=1/2$, by recalling the arguments used in \cite{GJ5}. As a consequence, the
fluctuations of the current are easily derived.

\subsection{Density fluctuations}

 Recall the definition of the density field given in \eqref{gamma density field} with $\gamma=1/2$. We want to obtain the limiting density field as a solution
of the stochastic partial differential equation \eqref{sB}.

For that purpose we introduce the martingales associated to $\mathcal{Y}_{t}^{n,\gamma}$ as in the previous section. The main difference on the limiting density field for
this strength of the asymmetry, comes from the limit of the integral term $\mathcal{A}_{t}^{n,\gamma}(H)$. We saw above that if $\gamma>1/2$ the somehow "stronger"
Boltzmann-Gibbs principle derived in \cite{G1}, tells us that $\mathcal{A}_{t}^{n,\gamma}(H)$ vanishes in $L^{2}(\mathbb{P}_{\nu_{\rho}})$ as $n\rightarrow{+\infty}$.

In the presence of a stronger asymmetry, a second order Boltzmann-Gibbs principle is needed and in \cite{GJ5} it was derived. This stronger replacement is derived through a multi-scale argument that was introduced in \cite{G1}. This multi-scale approach,
allows to obtain the desired replacement as a sequence of minor replacements in microscopic boxes that duplicate size at each step, until a
point in which the sum of the errors committed at each step is negligible and the last replacement holds at a box of small macroscopic size. In this macroscopic box, the non-trivial term can be identified as the square of $\mathcal{Y}_{t}^{n,\gamma}$. For details on this argument we refer the reader to \cite{GJ5}. The fundamental features of the
model that are used in order to derive this second order Boltzmann-Gibbs principle, are the sharp spectral gap bound for the dynamics restricted
to finite boxes, plus a second order expansion on the equivalence of ensembles. These results are quite general and this
argument can be applied for more general models than of exclusion type.

With this procedure, in \cite{GJ5} was shown that for any $H \in
\mc S(\bb R)$, $\mathcal{A}_{t}^{n,\gamma}(H)$ converges in $L^2(\mathbb{P}_{\nu_{\rho}})$ to
\begin{equation}\label{limit a}
\mathcal{A}_{t}(H)=\lim_{\epsilon\rightarrow{0}}-\frac{\beta''(\rho)}{2}\int_0^t\int_{\mathbb{R}} \mathcal{Y}_{s}(i_{\epsilon}(x))^2
H'(x) dx ds.
\end{equation}
 Here $i_{\epsilon}(x)(y)=\epsilon^{-1}1_{(x<y\leq{x+\epsilon})}$.

Applying the same arguments as above, in \cite{GJ5} it was shown that any limit point of $\mathcal{Y}_{t}^{n,\gamma}$ satisfies:
\begin{equation*}
\mathcal{Y}_{t}(H)-\mathcal{Y}_{0}^{\gamma}(H)=\int_{0}^{t}\mathcal{Y}_{s}\Big(\frac{\varphi_h'(\rho)}{2}\Delta H\Big)ds
+\mathcal{A}_{t}(H)+||\sqrt{\beta(\rho)}\nabla H||_{2} \mathcal{W}_{t}(H)
\end{equation*}
with $\mathcal{A}_{t}(H)$ given as in \eqref{limit a}. According to \cite{GJ5}, the limit density field $\mathcal{Y}_{t}$ is a weak solution of the KPZ equation
\eqref{sB}.

\subsection{Current fluctuations}

As mentioned above, having established the fluctuations of the density of particles one can obtain the fluctuations of the current of particles
$\mathcal{J}_{0}^{n,\gamma}(t)$. Following the route described above, it follows that
\begin{equation*}
\frac{1}{\sqrt{n'}}\Big\{\mathcal{J}_{0}^{n',\gamma}(t)-\mathbb{E}_{\nu_{\rho}}[\mathcal{J}_{0}^{n',\gamma}(t)]\Big\}
\end{equation*}
 converges to
\begin{equation*}
\mathcal{Y}_{t}(H_{0})-\mathcal{Y}_{0}(H_{0})
\end{equation*}
 where $\mathcal{Y}_{t}$ is the limit of the subsequence $\mathcal{Y}_{t}^{n',\gamma}$
and it is a weak solution of the KPZ equation \eqref{sB}.

So, the crossover regime for the current of particles occurs from Gaussian to a distribution which is given in terms of the solution of the KPZ
equation \eqref{sB}. As for the density of particles, the crossover of the current also occurs at strength asymmetry $an^{3/2}$.

In \cite{ACQ,SS2} the authors studied the weakly asymmetric simple exclusion process (which corresponds to taking $c(\cdot)\equiv1$ here) starting
from different initial conditions. There, the crossover regime for the current of particles is studied and the transition goes from Gaussian to
Tracy-Widom distribution. For the strength asymmetry $n^{3/2}$, \cite{SS2} identify the limit distribution of the current as the difference of two
Fredholm determinants which, from the random matrix theory, is known to converge to the Tracy-Widom distribution.

 In \cite{GJ5}, the stationary situation is studied and the limit of
the current is written in terms of the energy solution of the KPZ equation. The result is true for a general class of weakly asymmetric
exclusion processes. Here we show that the transition from the Edwards-Wilkinson class to the KPZ class is universal within this class of
processes and the strength asymmetry in order to have the crossover is precisely $an^{3/2}$. This is a step towards characterizing the
universality of the KPZ class.

\section*{Acknowledgements}
P.G. would like to thank the warm hospitality of Universit\'e Paris-Dauphine (France), where this work was initiated and to IMPA (Brazil) and Courant Institute of Mathematical Sciences (USA), where this work was finished. P.G. thanks to ``Funda\c c\~ao para a Ci\^encia e Tecnologia" for the research project PTDC/MAT/109844/2009: ``Non-Equilibrium Statistical Physics" and for the financial support provided by the Research Center of
Mathematics of the University of Minho through the FCT Pluriannual
Funding Program.

\end{document}